\documentclass[12pt]{amsart}

\usepackage{amsmath}
\usepackage{amsfonts}
\usepackage{amssymb}
\usepackage{amsthm}
\usepackage{hyperref}
\usepackage{enumerate}
\usepackage{cleveref}
\usepackage{mathrsfs}
\usepackage[top=0.9in, bottom=1.15in, left=1.1in, right=1.1in]{geometry}
 \usepackage{hyperref}
\hypersetup{colorlinks=true,linkcolor=blue,citecolor=magenta}
\newtheorem{theorem}{Theorem}

\newtheorem{proposition}[theorem]{Proposition}

\Crefname{conjecture}{Conjecture}{Conjectures}

\theoremstyle{plain}

\theoremstyle{plain}

\newcommand{\F}{\mathcal{F}}

\author{Robert Schneider} 
\address{Department of Mathematics and Computer Science\newline
Emory University\newline
400 Dowman Dr., W401\newline
Atlanta, Georgia 30322}
\email{robert.schneider@emory.edu}
\title[Fibonacci numbers and the golden ratio]{Fibonacci numbers and the golden ratio}
 
\begin{document}

\begin{abstract}
In this expository paper written to commemorate Fibonacci Day 2016, we discuss famous relations involving the Fibonacci sequence, the golden ratio, continued fractions and nested radicals, and show how these fit into a more general framework stemming from the quadratic formula.
\end{abstract}

\maketitle

\section{Fibonacci numbers}
%
%
The {\it Fibonacci numbers} are an interesting sequence of integers discovered by the prominent medieval mathematician Leonardo Fibonacci, related to the shapes of flower petals and tree branches, the birth rates of rabbits, and other natural phenomena, that shows up in many places in mathematics---you can even find the sequence in Pascal's triangle \cite{phi}. The $n$th Fibonacci number $F_n$ is defined for $n=0$ and $n=1$ by $$F_0=F_1=1,$$ and for $n\geq 2$ by 
\begin{equation}\label{fib}
F_n=F_{n-1}+F_{n-2}.
\end{equation}    
For example, we have 
\begin{flalign*}
F_2&=F_1+F_0=1+1=2\\F_3&=F_2+F_1=2+1=3\\F_4&=F_3+F_2=3+2=5,
\end{flalign*}
and so on. When you define a sequence of numbers like this, building up the $n$th term from the previous terms, it is called a {\it recursive sequence}, and an equation like (\ref{fib}) that produces the sequence is called a {\it recursion relation}. 

Recursion relations can lead to all sorts of lovely and intricate patterns, by plugging the terms back into themselves in creative ways. 
For example, the reader might like to try taking the definition (\ref{fib}) above, and repeatedly using the observation $$F_n=F_{n-1}+F_{n-2}=F_{n-1}+(F_{n-3}+F_{n-4})=F_{n-1}+F_{n-3}+(F_{n-5}+F_{n-6})$$
and so on, to show the following fact. 

\begin{proposition}
If $N$ is an odd number then 
$$F_N=F_0+F_2+F_4+F_6+...+F_{N-1}.$$
\end{proposition}
So ``$F_{odd}$'' can be written as the sum of all the smaller Fibonacci numbers of the form ``$F_{even}$''. For instance, $F_5=8$ can be written
$$F_5=F_0+F_2+F_4=1+2+5.$$

How about if $N$ is an {\it even} number above? The reader is encouraged to experiment with (\ref{fib}) to answer this question, and to discover other Fibonacci relations for yourself.
\  

\section{The golden ratio}
It is well known that the Fibonacci numbers are connected to another famous number, the {\it golden ratio} $\varphi$, which is related to the shapes of pineapples, seashells, and other objects in nature, as well as to fractals and other self-similar mathematical objects \cite{phi}. The golden ratio---
studied since the time of Euclid---is equal to the larger of the two roots\footnote{The quadratic formula says the roots, or ``zeros'', of $Ax^2+Bx+C$ are given by $\frac{-B\pm\sqrt{B^2-4AC}}{2A}$.} of the polynomial 
\begin{equation}\label{quadratic}
x^2-x-1.
\end{equation}
The quadratic formula gives the exact value\footnote{Interestingly, the other root $\frac{1-\sqrt{5}}{2}$ is equal to $-1/\varphi$.} 
$$\varphi=\frac{1+\sqrt{5}}{2}=1.6180... .$$ A beautiful fact is that $\varphi$ is related to its own reciprocal by a very simple formula.
\begin{proposition}\label{reciprocal}
The golden ratio is equal to its own reciprocal plus $1$:
\begin{equation*}
\varphi=1+\frac{1}{\varphi}
\end{equation*}
\end{proposition}
To see this, notice from expression (\ref{quadratic}) that $\varphi^2-\varphi-1=0$ implies $$\varphi^2=1+\varphi.$$ Now divide both sides of this equation by $\varphi$, to finish the proof.

Proposition \ref{reciprocal} gives us another kind of recursion relation. Observe that $\varphi$ is present 
on both the left-hand side, and in the denominator on the right. So we can plug the entire right-hand side of Proposition \ref{reciprocal} (which is equal to $\varphi$) back into itself in place of the $\varphi$ in the denominator:
$$\varphi=1+\frac{1}{1+\frac{1}{\varphi}}$$
There still is $\varphi$ in the denominator on the right-hand side, so we can repeat this substitution, and repeat it again, any number of times:
$$\varphi=1+\frac{1}{1+\frac{1}{1+\frac{1}{\varphi}}}=1+\frac{1}{1+\frac{1}{1+\frac{1}{1+\frac{1}{\varphi}}}}=1+\frac{1}{1+\frac{1}{1+\frac{1}{1+\frac{1}{1+\frac{1}{\varphi}}}}}=...$$
These stacked fractions-within-fractions are called {\it continued fractions}, and the numbers making up the fraction are called {\it coefficients}. Notice that at every stage, the left-hand side does not change: the continued fraction is still equal to $\varphi$. Proceeding in this manner forever, we get a very famous formula.
\begin{proposition}\label{phi}
We can write the golden ratio as an infinite continued fraction with all the coefficients equal to 1:
\begin{equation*}
\varphi=1+\frac{1}{1+\frac{1}{1+\frac{1}{1+\frac{1}{1+...}}}}
\end{equation*}
\end{proposition}
In the next section we will use this formula to show how $\varphi$ is related to the Fibonacci sequence. The interested reader is referred to Hardy and Wright \cite{HardyWright} for more about continued fractions.
\  

\section{Golden Fibonacci ratios}\label{Sect3}
The connection between the Fibonacci numbers $F_n$ and the golden ratio $\varphi$ is this.
\begin{proposition}
The ratio $\frac{F_n}{F_{n-1}}$ approaches 
$\varphi$ as $n$ increases.
\end{proposition}
From the very beginning of the Fibonacci sequence we see the ratio $F_n/F_{n-1}$ oscillates around $\varphi=1.6180...$, getting closer and closer to the golden ratio:
\begin{flalign*}
F_1/F_0&=1\\
F_2/F_1&=2\\
F_3/F_2&=3/2=1.5\\
F_4/F_3&=5/3=1.6666...\\
F_5/F_4&=8/5=1.6\  \  \text{(which is getting pretty close to $\varphi$ already)}
\end{flalign*}
Skipping up the sequence just a few terms, we have 
$$F_{10}/F_9=89/55=1.6181...,$$
which is indeed a very close approximation to $1.6180...$ (and it keeps getting better). 

It is easy to see why this is true, taking a recursive approach. Using the definition (\ref{fib}) of the $n$th Fibonacci number, 
we can rewrite this ratio as 
\begin{equation}\label{ratio}
\frac{F_n}{F_{n-1}}=\frac{F_{n-1}+F_{n-2}}{F_{n-1}}=1+\frac{F_{n-2}}{F_{n-1}}=1+\frac{1}{F_{n-1}/F_{n-2}}.
\end{equation}
By exactly the same principle, we can rewrite 
\begin{equation*}
\frac{F_{n-1}}{F_{n-2}}=1+\frac{1}{F_{n-2}/F_{n-3}},\  \  \  \frac{F_{n-2}}{F_{n-3}}=1+\frac{1}{F_{n-3}/F_{n-4}},
\end{equation*}
and so on, and substitute these one after another for $F_{n-1}/F_{n-2}$ in the right-hand side of equation (\ref{ratio}): 
$$\frac{F_{n}}{F_{n-1}}=1+\frac{1}{F_{n-1}/F_{n-2}}=1+\frac{1}{1+\frac{1}{F_{n-2}/F_{n-3}}}=1+\frac{1}{1+\frac{1}{1+\frac{1}{F_{n-3}/F_{n-4}}}}=...$$ 
Eventually we run out of Fibonacci numbers to put in the denominator, and end up with 
\begin{equation}\label{fibfrac}
{F_{n}}/{F_{n-1}}=1+\frac{1}{1+\frac{1}{1+\frac{1}{1+...+\frac{1}{F_1/F_0}}}}=1+\frac{1}{1+\frac{1}{1+\frac{1}{1+...+1}}},
\end{equation} 
where the final coefficient on the right is $1$ because $F_1/F_0=1$. Noting that the continued fraction on the right side gets longer and longer as $n$ increases (the number of coefficients increases proportionally), then as $n$ approaches $\infty$ the continued fraction in (\ref{fibfrac}) gets closer and closer to the infinite continued fraction in Proposition \ref{phi}, which equals $\varphi$.

That it turns out to be the exact number approached by ratios of consecutive Fibonacci numbers, which are interesting in their own right, is another impressive property of 
$\varphi$. 
\  

\section{Nested radicals} 
We want to point out one other beautiful (and exotic-looking) property of the golden ratio involving square roots, i.e., ``radicals''. The equation $\varphi^2-\varphi-1=0$ implies $\varphi^2=1+\varphi$; taking the square root of both sides leads to the following relation. 
\begin{proposition}\label{thmsqrt}
The golden ratio is equal to the square root of itself plus $1$: 
$$\varphi=\sqrt{1+\varphi}$$
\end{proposition} 
Viewing this as a recursion relation, just as with the continued fractions previously, we can now repeatedly substitute the entire right-hand side of Proposition \ref{thmsqrt} for the $\varphi$ under the radical on the right:
$$\varphi=\sqrt{1+\sqrt{1+\varphi}}=\sqrt{1+\sqrt{1+\sqrt{1+\varphi}}}=\sqrt{1+\sqrt{1+\sqrt{1+\sqrt{1+\varphi}}}}=...$$
These roots-within-roots are called {\it nested radicals}, the more obscure cousins of continued fractions. Continuing in this fashion forever, then $\varphi$ can be written as an infinite nested radical containing only $1$'s, like the continued fraction representation in Proposition \ref{phi}. 

So that we have the expressions all in one place, we collect this observation and our previous formulas for $\varphi$ in the following amazing statement.

\begin{proposition}\label{amazing}
The golden ratio is equal to
\begin{flalign*}
\varphi=\frac{1+\sqrt{5}}{2}=1+\frac{1}{1+\frac{1}{1+\frac{1}{1+\frac{1}{1+...}}}}=\sqrt{1+\sqrt{1+\sqrt{1+\sqrt{1+...}}}}.
\end{flalign*}
\end{proposition} 

In the next section, we will write a very similar system of equations for the roots of {\it any} quadratic polynomial with rational coefficients. 
\\   

\section{The quadratic formula and beyond}
Here we will see there is a more general framework containing the previous results about the golden ratio, Fibonacci numbers, continued fractions and nested radicals---all going back to the quadratic formula. 

Now, the $\pm$ sign in the quadratic formula yields two roots, a ``plus'' root and a ``minus'' root. For rational numbers $a$ and $b$, let $\varphi(a,b)$ 
denote the ``plus'' root of the polynomial \begin{equation}\label{polyn2}
x^2-ax-b.
\end{equation}
Then the quadratic formula 
gives the exact value\footnote{We can use the ``minus'' root $\frac{a - \sqrt{a^2+4b}}{2}=a-\varphi(a,b)=-b/\varphi(a,b)$ for similar results to those that follow, but with some $\pm$ sign changes. The reader is encouraged to work this case out, too.} for this number: 
\begin{equation*}
\varphi(a,b)=\frac{a + \sqrt{a^2+4b}}{2}
\end{equation*}
The golden ratio is the special case $\varphi=\varphi(1,1)$. 
In fact, these numbers $\varphi(a,b)$ possess many of the nice properties enjoyed by $\varphi$. 
For instance, it follows from (\ref{polyn2}) that 
$$\varphi(a,b)^2-a\varphi(a,b)-b=0.$$ Then similar steps to those we applied to $\varphi$ yield a pair of familiar-looking equalities.

\begin{proposition}
We can write $\varphi(a,b)$ in the following ways:
\begin{flalign*}
\varphi(a,b)&=a+\frac{b}{\varphi(a,b)}\  \  \  \text{(so long as $\varphi(a,b)\neq 0$ in the denominator)}\\
\\
&=\sqrt{b+a\varphi(a,b)}
\end{flalign*}
\end{proposition}

Notice how the case $a=b=1$ reduces to Propositions \ref{reciprocal} and \ref{thmsqrt}. Using these two identities, and following the exact steps that proved the corresponding equations for $\varphi$ previously, we can generalize Proposition \ref{amazing} as follows. We omit the proofs, however, as the reader might enjoy working out the details yourself. 

\begin{proposition}\label{amazing2}
The number $\varphi(a,b)$ is equal to
\begin{flalign*}
\varphi(a,b)=\frac{a + \sqrt{a^2+4b}}{2}=a+\frac{b}{a+\frac{b}{a+\frac{b}{a+\frac{b}{a+...}}}}=\sqrt{b+a\sqrt{b+a\sqrt{b+a\sqrt{b+...}}}}.
\end{flalign*}
\end{proposition} 

We can also define a generalization of the Fibonacci sequence that connects with $\varphi(a,b)$. Let us define $F_n(a,b)$ by $$F_0(a,b)=1,\  F_1(a,b)=a,$$ and for $n\geq 2$ by 
\begin{equation}\label{faux}
F_n(a,b)=aF_{n-1}(a,b)+bF_{n-2}(a,b).
\end{equation}
The usual Fibonacci numbers are the special case $F_n=F_n(1,1)$. It turns out these $F_n(a,b)$ behave similarly to $F_n$, as we see in the following statement.

\begin{proposition}\label{asymp2}
The ratio $\frac{F_n(a,b)}{F_{n-1}(a,b)}$ approaches $\varphi(a,b)$ as $n$ increases.
\end{proposition}
 
This proposition follows from the definition (\ref{faux}) with the observation 
\begin{equation}\label{genfib}
\frac{F_n(a,b)}{F_{n-1}(a,b)}=\frac{aF_{n-1}(a,b)+bF_{n-2}(a,b)}{F_{n-1}(a,b)}=a+\frac{b}{F_{n-1}(a,b)/F_{n-2}(a,b)},
\end{equation}
exactly like the proof of Proposition \ref{fibfrac}, so that after repeatedly using (\ref{genfib}) we end up with 
$$\frac{F_n(a,b)}{F_{n-1}(a,b)}=a+\frac{b}{a+\frac{b}{a+\frac{b}{a+...+\frac{b}{F_1(a,b)/F_0(a,b)}}}}=a+\frac{b}{a+\frac{b}{a+\frac{b}{a+...+\frac{b}{a}}}}.$$
As $n$ increases, the far right side of this equation looks more and more like the infinite continued fraction in Proposition \ref{amazing2}, that is, more and more like $\varphi(a,b)$. 

What if we define $F_0(a,b)$ to be a different number from $1$? 
What if we change $F_1(a,b)$ too? Does this affect Proposition \ref{asymp2}? And what if we use more than two terms in the recursion relation (\ref{faux}), for instance, if we define $F_n(a_1,a_2,a_3,...,a_k)$ using $k$ terms? Is there a special number $\varphi(a_1,a_2,a_3,...,a_k)$ 
associated with such a sequence, 
as $\varphi(a,b)$ is with $F_n(a,b)$? 
The reader is encouraged to 
experiment with 
Fibonacci-like sequences, and also to check out important variations on the Fibonacci numbers, such as Lucas numbers (see \cite{HardyWright}), that produce other interesting relations.

%

\end{document}